\documentclass[11pt,reqno]{amsart}
\usepackage{amsfonts}
\usepackage{amsmath}
\usepackage{amssymb}
\usepackage{amscd}
\usepackage[latin2]{inputenc}
\usepackage[T1]{fontenc}

\newcommand{\intfrac}[2]{\genfrac{\lfloor}{\rfloor}{}{1}{#1}{#2}}

\newtheorem*{theorem}{Main Theorem}

\numberwithin{equation}{section}

\begin{document}
\title[Number of singular points]{Number of singular points of an annulus in 
$\mathbb{C}^2$}
\author{Maciej Borodzik}
\address{Institute of Mathematics, University of Warsaw, ul. Banacha 2,
02-097 Warsaw, Poland}
\email{mcboro@mimuw.edu.pl}
\author{Henryk \.Zoł{}ądek}
\address{Institute of Mathematics, University of Warsaw, ul. Banacha 2,
02-097 Warsaw, Poland}
\email{zoladek@mimuw.edu.pl}
\date{April 15, 2010}
\subjclass{primary: 14H50, secondary: 14R10, 14B05}
\keywords{annulus, cuspidal singular point, codimension}
\thanks{The first author is supported by Foundation for Polish Science and
both authors are supported by Polish MNiSzW Grant No 2 P03A 015 29}

\begin{abstract}
Using BMY inequality and a Milnor number bound we prove that any algebraic
annulus $\mathbb{C}^*$ in $\mathbb{C}^2$ with no self--intersections can
have at most three cuspidal singularities.
\end{abstract}

\maketitle

\section{Introduction}

The problem of classification of curves in $\mathbb{C}^{2}$ of fixed
topological type up to an algebraic automorphism of $\mathbb{C}^{2}$ is in
general very difficult. One of the most important result in this domain is
the Abhyankar--Moh--Suzuki theorem ([AM, Suz]) stating that any algebraic
curve in $\mathbb{C}^{2}$ that is diffeomorphic to a disk is in fact
algebraically isomorphic to a line. Another one, due to M. Zaidenberg and V.
Lin [LZ1], says that any curve homeomorphic to a disk is algebraically
equivalent to a curve of the type $x^{p}=y^{q}$ for $p,q$ coprime.

In [BZ1, BZ2] we developed an efficient method in some other particular
cases: namely we studied rational curves with one place at infinity and one
double point (topological immersions of $\mathbb{C}$ in $\mathbb{C}^{2}$
with one finite self-intersection) in [BZ1] and annuli (topological
embeddings of $\mathbb{C}^{\ast }$ in $\mathbb{C}^{2}$) in [BZ2]. A list of
44 possible cases was found and it was claimed that the list is complete.
The claim boils down to the validity of certain conjecture, strongly related
to the unobstructedness problem of [FZ]. (We also refer the reader to the
paper [CKR] where a partial classification of annuli is given.)

It turns out that our method, even without assuming the above--mentioned
conjecture, can be applied to prove a conjecture by Lin and Zaidenberg [LZ2]
specified to annuli. The latter conjecture states that any algebraic curve
in $\mathbb{C}^{2}$ with the first Betti number equal to $r$ can have at
most $2r+1$ singular points. In the present paper we prove the following
theorem, which confirms the Zaidenberg--Lin conjecture for annuli.

\begin{theorem}
Any algebraic curve in $\mathbb{C}^2$ homeomorphic to $\mathbb{C}^*$ has at
most three singular points.
\end{theorem}

The method of the proof is as follows. We use a notion of codimension of a
singular point (see [BZ3]). This is the number of conditions for a
parametric curve required so that this curve has a given singularity (up to
a topological equivalence). A parameter count argument would give the bound
for the sum of codimensions over all singular points of the given curve by
the dimension of the space of parametric curves. This dimension depends
linearly on the degree of the curve under consideration. In [BZ2, Conjecture
3.7] we conjectured such bounds. While we do not have the proof of these
bounds, we noted that a slightly weaker codimension bound can be obtained
using Bogomolov--Miyaoka--Yau (BMY) inequality (compare [BZ3]). This bound
being insufficient to prove that the list in [BZ2] is complete, at least
without an additional work, yet is suitable to verify that an annulus cannot
have more than three finite singular points.

We believe that our methods can settle the conjecture for all rational
curves in $\mathbb{C}^{2}$. However the computations in the general case
seem to be highly complex. In the case of affine plane curves of arbitrary
genus with one place at infinity some estimates for the number of singular
points have been recently obtained in [Bor]. \medskip

\textbf{Acknowledgments.} This paper was motivated and partially written
during the workshop \textquotedblleft Affine algebraic
geometry\textquotedblright\ in Oberwolfach. We are grateful to the
organizers for inviting us to this conference. We would like to thank M.
Koras and V. Lin for stimulating discussions. The first author thanks P.
Russell for interesting discussions and an invitation to visit the McGill
University.

\bigskip

\section{Invariants of singular points}

Here we present some notions and estimates from [BZ1, BZ2, BZ3].\medskip

\subsection{Local invariants of singularities of curves}

Let $(A,0),$\linebreak\ $A=\left\{ f(x,y)=0\right\} \subset \mathbb{C}^{2},$
be a germ of a reduced plane curve near its singular point.

The first invariants of this singularity are: the \textit{number of branches}
(irreducible components), denoted by $k,$ and the \textit{multiplicity},
denoted by \textrm{mult}$_{0}A$\linebreak $=m$. The latter is the order of
the first nonzero term in the Taylor expansion of the defining function $f$.
In this work we consider only the cases with $k=1$ (cuspidal singularities)
and $k=2$ (for an annulus it may occur at infinity).

Next invariant is the \textit{external codimension} of a singularity,
denoted $ext\nu $ and defined as follows.

Let 
\begin{equation}
x=\tau ^{m},\text{ \ \ }y=c_{1}\tau +c_{2}\tau ^{2}+\ldots   \label{2.1}
\end{equation}%
be the Puiseux expansion of $A$ in the cuspidal case ($k=1$). In the space
of germs as above (i.e. with fixed multiplicity $m)$\ strata of topological
equivalence (or so-called $\mu =$\textrm{const} strata) are defined by
vanishing of some number of certain Puiseux quantities $c_{j}$ and by
nonzero some other Puiseux quantities; in [BZ2, Section 2.II] and [BZ3,
Section 2.1] the quantities $c_{j}$ which appear in descriptions of these
strata are called the \textit{essential Puiseux quantities}. The number of
vanishing essential Puiseux quantities is the $y-$\textit{codimension}
denoted by $\nu $ (see [BZ3, Section 2.1]). This can be explained in terms
of the so-called \textit{topologically arranged Puiseux expansion} 
\begin{equation*}
y=\left( d_{0}x^{\frac{n_{0}}{1}}+\ldots \right) +\left( d_{1}x^{\frac{n_{1}%
}{m_{1}}}+\ldots \right) +\ldots +\left( d_{r}x^{\frac{n_{r}}{m_{1}\ldots
m_{r}}}+\ldots \right) ,
\end{equation*}%
where $n_{j}$ and $m_{j}$ are positive integers (with $1=m_{0}$ and $m_{j}>1$
for $j\geq 1$) such that $\gcd \left( m_{j},n_{j}\right) =1$, $m_{1}\ldots
m_{r}=m$, the nonzero coefficients $d_{1},\ldots ,d_{r}$ constitute a part
of the essential Puiseux quantities and the dots in the $j^{\mathrm{th}}$
summand mean terms with $x^{\frac{i}{m_{0}\ldots m_{j}}}.$ Here the first
(inessential) summand can be absent and the pairs $\left( m_{j},n_{j}\right)
,$ $j\geq 1$, are known as the \textit{Puiseux pairs}. The other essential
Puiseux quantities, i.e.other than $d_{j}=c_{n_{j}m_{j+1}\ldots m_{r}}$ for $%
j=1,\ldots ,r,$ correspond to those terms $c_{i}x^{i/m}$ whose potential
presence would change the essential part of the above topologically arranged
Puiseux expansion.

For example, in the case $m=2$\ the strata of topological equivalence are
defined by $c_{1}=c_{3}=\ldots =c_{2\nu -1}=0\not=c_{2\nu +1}$\ in (2.1). In
the case $m=4$\ the conditions $c_{1}=c_{2}=c_{3}=0\not=c_{5},$ $%
c_{1}=c_{2}=c_{3}=c_{5}=c_{7}=c_{9}=0\not=c_{6}c_{11}$ and $%
c_{1}=c_{2}=c_{3}=c_{5}=c_{6}=c_{7}=c_{9}=0\not=c_{10}c_{11}$ define three $%
\mu =$\textrm{const} strata with $\nu =3,$ $\nu =6$ and $\nu =7$
respectively. We see that the name `essential' for a Puiseux quantity
sometimes depends on the stratum (like for $c_{10}$ above), but the
quantities $c_{m},c_{2m},\ldots $\ are always\ inessential.\ 

In the cuspidal case we put 
\begin{equation*}
ext\nu =\nu +m-2;
\end{equation*}%
the additional contribution to $ext\nu $ arises from the conditions for the $%
m-1$\ first derivatives of $x$ (with respect to a parameter $t$ on the
curve) to vanish at the singular point.

(Formally one obtains $m-1$ independent conditions, i.e. in the space of
local parametric curves $\left( \mathbb{C},0\right) \longmapsto \left( 
\mathbb{C}^{2},0\right) $. However, in the space of global parametric
curves, like in (2.3) below, the positions of their singular points are not
fixed. So the condition $dx/dt=0$ is just the equation for values of the
parameter at the singular points and, as such, it does not enter into the
collection of `external conditions' for the singularity.)

In the two branches case, $A=A_{1}+A_{2}$ with the multiplicities $m_{1}$
and $m_{2}$, besides the $y-$codimensions $\nu (A_{1})$ and $\nu (A_{2}),$\
we have also the \textit{tangency codimension} $\nu _{\tan }$\ between the
two branches. It is the number of inessential Puiseux quantities and nonzero
essential Puiseux quantities in the common part of the Puiseux expansions of
the two branches (we choose the roots of unity of orders $m_{1}$ and $m_{2}$
to make this common part as long as possible). $\nu _{\tan }$\ is a
topological invariant of the singularity, because it controls the
intersection index of the branches.

For example, if the Puiseux expansions of the two branches are $y=\alpha
x^{3/2}+x^{2}+\ldots $\ and $y=\beta x^{3/2}-x^{2}+\ldots $\ (with $%
m_{1}=m_{2}=2$) then $\nu (A_{1})=\nu (A_{2})=1$\ and $\nu _{\tan }=1$\ when 
$\alpha ^{2}\not=\beta ^{2}$ and $\nu _{\tan }=2$ otherwise.

Here we put 
\begin{equation*}
ext\nu =\nu (A_{1})+\nu (A_{2})+\nu _{\tan }+m_{1}+m_{2}-2.
\end{equation*}

\medskip

There exists another interpretation of the external codimension. Namely, we
take the \textit{minimal normal crossing resolution} of the singular point $%
\pi :\left( V,D\right) \rightarrow \left( U,A\right) ,$ where $U$ is a
neighborhood of the origin in $\mathbb{C}^{2}.$ Letting $E=E_{1}+\ldots
+E_{l}$ be the \textit{exceptional divisor} with components $E_{j},$ we
consider the vector space 
\begin{equation*}
\mathrm{Vect}(E)=\mathbb{Q}E_{1}\oplus \ldots \oplus \mathbb{Q}E_{l}
\end{equation*}%
equipped with the intersection form. Then the strict transform $\widetilde{A}
$ of $A,$ as well as $D,$ the reduced total transform of $A$, are
interpreted as elements of $\mathrm{Vect}(E).$ We have also the \textit{%
local canonical divisor} $K$ defined by the relations $E_{j}(K+E_{j})=-2$
(see [OZ2]).

The following result was proved by S. Orevkov [Or] in the cuspidal case and
in [BZ3, Proposition 4.1] in general. Orevkov calls the quantity $K(K+D)$
the \textit{rough M--number}.

\medskip

\textbf{Proposition 2.1.} \textit{We have} 
\begin{equation*}
ext\nu =K(K+D).
\end{equation*}

\medskip

A classical invariant of singularity is the \textit{number of double points}%
, denoted by $\delta $ (sometimes called the delta invariant). In the
cuspidal case it equals $\mu /2,$ where $\mu $\ is the \textit{Milnor number}
of the singularity. Generally it is the number of double points of a
parametric deformation of the curve $A:$\ we take a\ map from a disjoint
union $\coprod \left\{ \left\vert z\right\vert <\varepsilon \right\} $\ of
complex discs\ to $\mathbb{C}^{2}$\ which is a small generic perturbation of
the normalization map. In this sense we can interpret $\delta $ as the
number of double points which are hidden at the singularity.

For example, for the $\mathbf{A}_{\mu }$ singularity $y^{2}=x^{\mu +1}$ we
have $\delta =\mu /2$\ if $\mu $\ is even and $\delta =(\mu +1)/2$\ if $\mu $%
\ is odd.

The following inequality was proved in [BZ1, Proposition 2.9 and Proposition
2.16].

\medskip

\textbf{Proposition 2.2.} \textit{If the number of branches is }$k=1$\textit{%
\ or }$k=2$\textit{\ and }$m=\mathrm{mult}_{0}A$\textit{\ then }%
\begin{equation*}
2\delta \leq m(ext\nu -m+k+1).
\end{equation*}

\medskip

In the above vector space $\mathrm{Vect}(E)$, related with the resolution of
singularity, we can use the local \textit{Zariski--Fujita decomposition}
[OZ2]%
\begin{equation*}
K+D=P+N,
\end{equation*}%
where $P$ is the positive and $N$ the negative part of $K+D$ (with respect
to the intersection form). Then we define the \textit{excess} of the
singular point as 
\begin{equation}
\eta :=-N^{2}\geq 0.  \label{2.2}
\end{equation}%
This is also a topological invariant, because it is defined via the
intersection form on the space $\mathrm{Vect}(E)$.

The following result follows from a rather subtle analysis of the
intersection form via dual graph by Orevkov and Zaidenberg [OZ2] (see also
[BZ3, Proposition 4.2]). Below we use the notations $\left\lfloor
x\right\rfloor =\max \{n:n\in \mathbb{Z},$ $n\leq x\},$ $\lceil x\rceil
=\min \{n:n\in \mathbb{Z},$ $x\leq n\}$. \medskip

\textbf{Proposition 2.3.} \textit{If }$(m,n)$\textit{\ is the first
characteristic pair of an unibranched singularity then its excess (2.2)
satisfies }%
\begin{equation*}
\eta \geq \left( \lceil m/n\rceil -m/n\right) +\left( \lceil n/m\rceil
-n/m\right) .
\end{equation*}%
\medskip

\subsection{Invariants of the annuli}

Consider an annulus $C$ given in parametric form by 
\begin{equation}
\begin{split}
x& =\varphi (t)=t^{p}+a_{1}t^{p-1}+\dots +a_{p+r}t^{-r} \\
y& =\psi (t)=t^{q}+b_{1}t^{q-1}+\dots +b_{q+s}t^{-s},
\end{split}
\label{2.3}
\end{equation}%
where $a_{p+r}b_{q+s}\neq 0$. The numbers $p$, $q$, $r$ and $s$ are integers
and we can assume that $p,s>0$ (since we have a topological embedding of $%
\mathbb{C}^{\ast }).$

Such an annulus may have several finite singular points corresponding to the
values $t_{1}$, $\dots $, $t_{N}$ of $t$. They are all cuspidal. The above
invariants associated with each point $t_{i}$ are denoted by $m_{i},$ $%
ext\nu _{i}$, $\delta _{i}$ and $\eta _{i}$.

We denote by $\nu _{\infty }$ the so-called \textit{subtle codimension} of
the branch of $C$ as $t$ goes to infinity, which is the codimension of the
topological equivalence stratum in the space of germs of the form $x=\tau
^{-p},$ $y=\tau ^{-q}+c_{1}\tau ^{-q+1}+\ldots ,$\ $\tau \rightarrow 0$\
(compare [BZ3, Definition 2.6]). Analogously we define the \textit{subtle
codimension} $\nu _{0}$ of the branch of $C$ as $t\rightarrow 0$.

The last invariant of the curve $C$ is the \textit{tangency codimension} $%
\nu _{\tan }$ at infinity. More precisely, if $ps\neq rq$ then the two
branches of $C$ do not intersect and we put $\nu _{\tan }=0.$ If $ps=rq$
then $\nu _{\tan }$ is defined as above for a two branches singularity. We
use the notion of $\nu _{\tan }$ only in Section 3.3 (Case B2). Sometimes we
will use the notation%
\begin{equation}
\nu _{\inf }=\nu _{0}+\nu _{\infty }+\nu _{\tan }.  \label{2.4}
\end{equation}

For the purpose of proving the Main Theorem, the above quantities are not
that important as the inequalities that relate them. The first identity,
which is a direct consequence of the standard genus formula (or the Poincaré%
--Hopf formula), can be found in [BZ2, Proposition 2.9 and Eq. (2.11)].

\medskip

\textbf{Proposition 2.4.} \textit{A generic curve of the form (2.3) has} 
\begin{equation*}
2\delta _{\max }:=(p+r-1)(q+s-1)+|ps-rq|-p^{\prime }-r^{\prime }+1
\end{equation*}%
\textit{finite simple double points, where}%
\begin{equation}
p^{\prime }=\gcd (p,q),\text{ \ \ }r^{\prime }=\gcd (r,s).  \label{2.5}
\end{equation}

\medskip

Since we are interested in the annuli, which by definition do not have
self-intersections, the $\delta _{\max }$ double points must be hidden at
singular points and/or at infinity: 
\begin{equation*}
\delta _{\max }=\sum \delta _{i}+\delta _{\inf }.
\end{equation*}

The numbers $\delta _{i}$ are estimated directly in Proposition 2.2,%
\begin{equation*}
2\delta _{i}\leq m_{i}(ext\nu _{i}-m_{i}+2).
\end{equation*}%
From that proposition we find also a bound for the number of double points
hidden at infinity (see [BZ2, Proposition 2.29]): 
\begin{eqnarray*}
2\delta _{\inf } &\leq &p^{\prime }\nu _{\infty }+r^{\prime }\nu _{0}\text{
\ if \ }ps\not=qr. \\
2\delta _{\inf } &\leq &(p^{\prime }+r^{\prime })(\nu _{\inf }+1)\text{ \ if
\ }ps=qr.
\end{eqnarray*}

We introduce the following quantity:%
\begin{equation}
\begin{array}{ll}
\mathcal{E}=\sum_{i=1}^{N}m_{i}(ext\nu _{i}-m_{i}+2)+p^{\prime }\nu _{\infty
}+r^{\prime }\nu _{0} & \text{ if }ps\not=qr, \\ 
\mathcal{E}=\sum_{i=1}^{N}m_{i}(ext\nu _{i}-m_{i}+2)+(p^{\prime }+r^{\prime
})(\nu _{\inf }+1) & \text{ if }ps=qr.%
\end{array}
\label{2.6}
\end{equation}%
By the above local estimates the inequality 
\begin{equation}
\Delta :=2\delta _{\max }-\mathcal{E}\leq 0  \label{2.7}
\end{equation}%
holds for an annulus of the form (2.3). The quantity $\Delta $ is called the 
\textit{reserve} in [BZ2, Section 2.1].

\medskip

Next we would like to bound the sum of codimensions. The bound depends on
values of the exponents $p$, $q$, $r$ and $s$. \medskip

\textbf{Definition 2.1.} A curve $C$\ given in (2.3) is of

\qquad --- \textit{type} $\binom{+}{+}$ if $0<p<q$ and $0<r<s$, $p+r<q+s$;

\qquad --- \textit{type} $\binom{-+}{+-}$ if $0<q<p$ and $0<r<s$, $p+r\leq
q+s$;

\qquad --- \textit{type} $\binom{-}{+}$ if $r<0$ and $q>0$;

\qquad --- \textit{type} $\binom{-}{-}$ if $r<0$ and $q<0$, $p+r\leq q+s$.

\medskip

Recall that with the open surface $V_{0}=\mathbb{C}^{2}\setminus C$ we can
associate its \textit{logarithmic Kodaira dimension} $\bar{\kappa}(\mathbb{C}%
^{2}\setminus C)$. It is defined via the normal crossing completion $F$ of $%
V_{0}$\ such that $V=V_{0}\cup F$ is smooth projective surface. Then 
\begin{equation*}
\bar{\kappa}(V_{0})=\lim \sup \log h^{0}(V,n(K_{V}+F))/\log n.
\end{equation*}%
If $\bar{\kappa}(V_{0})=2$ then we say that the surface $V_{0}$ is of 
\textit{general type}. I. Wakabayashi [Wa] calculated the logarithmic
Kodaira dimension of $\mathbb{C}^{2}\setminus C$\ in some important cases.
From [Wa] one can deduce, in particular, the following fact.

\medskip

\textbf{Proposition 2.5.} \textit{If an annulus }$C$\textit{\ has more than
three finite singular pints then the surface }$\mathbb{C}^{2}\setminus C$%
\textit{\ is of general type}.

\medskip

The codimension bounds we give below were proved in [BZ3, Theorem 4.3]; they
essentially rely upon the Bogomolov--Miyaoka--Yau inequality (which was also
used by Zaidenberg and Orevkov [Or, OZ2]). Here we state only the result.

Introduce the quantity%
\begin{equation}
S:=\sum_{i=1}^{N}(ext\nu _{i}+\eta _{i})+\nu _{\inf },  \label{2.8}
\end{equation}%
which can be regarded roughly as the sum of local codimensions.

\medskip

\textbf{Proposition 2.6.} \textit{Let }$C$\textit{\ be an annulus given by
(2.3) homeomorphic to }$\mathbb{C}^{\ast }$\textit{\ and such that its
complement in }$\mathbb{C}^{2}$\textit{\ is of general type. Then depending
on the type of the annulus we have:}

\textit{(a) for type }$\binom{+}{+}$%
\begin{equation}
S\leq p+r+q+s+1-\min (\lfloor q/p\rfloor ,\lfloor s/r\rfloor )\leq p+r+q+s%
\mathit{;}  \label{2.9}
\end{equation}

\textit{(b) for type }$\binom{-+}{+-}$%
\begin{equation}
S\leq p+r+q+s+1\mathit{;}  \label{2.10}
\end{equation}

\textit{(c) for type }$\binom{-}{+}$%
\begin{equation}
S\leq p-|r|+q+s+2+\lfloor \left( |r|-1\right) /s\rfloor -\lfloor q/p\rfloor 
\mathit{;}  \label{2.11}
\end{equation}

\textit{(d) for type }$\binom{-}{-}$%
\begin{equation}
S\leq p-|r|-|q|+s+3+\lfloor \left( |r|-1\right) /s\rfloor +\lfloor \left(
|q|-1\right) /p\rfloor .  \label{2.12}
\end{equation}

\medskip

For the multiplicities $m_{i}$ and the excesses $\eta _{i}$ (see (2.2)) of
singular points we have the following bounds.

$\medskip $

\textbf{Lemma 2.1.} 
\begin{equation*}
\sum_{i=1}^{N}(m_{i}-1)\leq \min (p+r,q+s).
\end{equation*}%
\textit{In particular, }$N\leq p+r$.

\medskip

\textit{Proof}. Assume that $p+r\leq q+s$. Then $\dot{x}=d\varphi
/dt=R(t)t^{-r-1}$, where $R(t)$ is a polynomial of degree $p+r$ (see (2.3)).
If $n_{i}-1$ is the order of $d\varphi /dt$ at the $i^{\mathrm{th}}$
singular point, then clearly $\sum (n_{i}-1)\leq p+r$ and $m_{i}\leq n_{i}$.
The second statement is obvious. $\qed$

$\medskip $

\textbf{Lemma 2.2.} \textit{(a) }$\eta _{i}>1/2,$ \textit{thus if }$N\geq 4$%
\textit{\ then }$\sum \eta _{i}>2.$\textit{\ }

\qquad \textit{(b) If the multiplicity of a singular point is }$m_{i}=2$%
\textit{\ then }$\eta _{i}\geq 5/6.$

\qquad \textit{(c) If }$N\geq 4$\textit{\ and }$\min (p+r,q+s)\leq 5$\textit{%
\ then }$\sum \eta _{i}>3.$

$\medskip $

\textit{Proof}. The first two assertions follow directly from Proposition
2.3. In (c) we must have $m_{i}=2$ for at least three singular points if
their number $N=4$ and all $m_{i}=2$ if $N=5.$ Therefore, by Proposition
2.3, either $\sum \eta _{i}>3\cdot \frac{5}{6}+\frac{1}{2}=3$ or $\sum \eta
_{i}\geq 5\cdot \frac{5}{6}>4.\ \qed$

\medskip

Two technical statements below turned out useful. The first one is often
used in [BZ2, Lemma 5.3].

\medskip

\textbf{Lemma 2.3.} \textit{If }$ps-rq\neq 0$\textit{\ then the quantity }%
\begin{equation*}
\det \,^{\prime }:=|ps-rq|-p^{\prime }-r^{\prime }+1
\end{equation*}%
\textit{is a non-negative integer}.

\medskip

The second lemma gives a partial answer to the problem of finding the best
parametrization of an annulus given by (2.3). In fact if, say, $%
x=t^{2}+\dots +t^{-6}$ and $y=t^{4}+\dots +t^{-9}$ we can ask whether it is
reasonable to apply a de Jonqui\`{e}re transform $y\rightarrow y-x^{2}$ to
reduce the order of $y$ at $t\rightarrow \infty $ at the cost of increasing
its order as $t\rightarrow 0$. We prove that there exists (maybe not unique)
way of choosing an automorphism of $\mathbb{C}^{2}$ that suits best to our
estimates.

\medskip

\textbf{Definition 2.2.} A curve $C$\ is called \textit{ugly} if one of the
following holds:

\qquad --- it is of type $\binom{+}{+}$ , $q/p\in \mathbb{Z}$ and $r<p$;

\qquad --- it is of type $\binom{-+}{+-}$ and either $p/q\in \mathbb{Z}$ and 
$s<q$ or $s/r\in \mathbb{Z}$ and $p<r$;

\qquad --- it is of type $\binom{-}{+}$ , $p/q\in \mathbb{Z}$ and $s<q$.

Otherwise the curve $C$ is called \textit{handsome}.

\medskip

\textbf{Lemma 2.4.} \textit{Any curve as in (2.3) can be transformed to a
handsome one by applying a Cremona automorphism of }$\mathbb{C}^{2}$\textit{%
\ and, possibly, the change} $t\rightarrow 1/t$.

\medskip

A straightforward proof is presented in [BZ2, Proposition 2.45].

\medskip

\subsection{Scheme of the proof of Main Theorem}

We can order the singular points of $C$ so that $m_{1}\geq m_{2}\dots \geq
m_{N}$.

Recall that we must rule out the possibility $N\geq 4$. But one quickly
realizes that considering the case $N=4$ is sufficient. As in [BZ2] the
estimates become easier when $N$ grows. For example, the codimension bound
is stronger already for $N=5$.

We split the proof into following five cases:

\begin{itemize}
\item[A] Type $\binom{+}{+}$ with $ps\neq rq$,

\item[B] Type $\binom{+}{+}$ with $ps=rq$,

\item[C] Type $\binom{-+}{+-}$,

\item[D] Type $\binom{-}{+}$,

\item[E] Type $\binom{-}{-}$.
\end{itemize}

Each case can be split in turn into two subcases.

\begin{itemize}
\item[$(1)$] We assume that double points hide at finite singular points.
This means that the quantity $\mathcal{E}$ from (2.6), which we try to
maximize, is greatest when $\delta _{\inf }=0.$ Then it is easy to see that $%
\mathcal{E}$ is maximal possible, when the multiplicity $m_{1}$ and the
external codimension $ext\nu _{1}$ are maximal, and the other multiplicities
and $ext\nu $ numbers, including $\nu _{\inf }=\nu _{0}+\nu _{\infty }+\nu
_{\tan }$ (see (2.4)) are minimal. Here we have $m_{1}\geq \max (p^{\prime
},r^{\prime })$ (see (2.5)) in cases A, C, D, E and $m_{1}\geq p^{\prime
}+r^{\prime }$ in case B.

\item[$(2)$] We assume that $\nu _{\inf }$ is large, so double points hide
at infinity (i.e. $\sum \delta _{j}$ is small relatively to $\delta _{\inf }$%
). Then $\mathcal{E}$ is maximal if all codimensions of singularities at
finite distance are minimal and the codimension at infinity is maximal
possible. Here either $p^{\prime }$\ or $r^{\prime }$ exceeds $m_{1}$\ in
cases A, C, D, E or $p^{\prime }+r^{\prime }>m_{1}$\ in case B.\ 
\end{itemize}

In all cases we shall strive to prove that the reserve $\Delta >0$, which
contradicts inequality (2.7). To simplify arguments we will assume that the
curve has precisely $N=4$ singular points.

\bigskip

\section{Proof of Main Theorem}

\subsection{Cases A1 and B1.}

These two cases are very similar. Here $\mathcal{E}$ is maximal if $%
m_{2}=m_{3}=m_{4}=2$, $ext\nu _{2}=ext\nu _{3}=ext\nu _{4}=1$ and $\nu
_{\infty }=\nu _{0}=0$. Hence $m_{1}\leq p+r-2$, and $ext\nu _{1}\leq
p+r+q+s-6$ ($-3$ coming from $ext\nu _{2}+ext\nu _{3}+ext\nu _{4}$ and
another $-3$ from $\sum \eta _{i}>2$, see (2.8)--(2.9) and Lemma 2.2 (a)).
Therefore, by (2.6), 
\begin{equation*}
\mathcal{E}\leq (p+r-2)(q+s-2)+6.
\end{equation*}

Let $ps\neq rq$. We have $2\delta _{\max }\geq (p+r-1)(q+s-1)$ (see
Proposition 2.4 and Lemma 2.3). Thus 
\begin{equation*}
\Delta =2\delta _{\max }-\mathcal{E}=(p+r+q+s)-9.
\end{equation*}%
But $p+r\geq 4$ by Lemma 2.1. Moreover, $q>p$ and $s>r$ by Definition 2.1.
Hence $q+s\geq 6$ and $\Delta >0$, so that there is no such curve $C$ in
this case (compare (2.7)).

Let $ps=rq$, so $2\delta _{\max }=(p+r-1)(q+s-1)-(p^{\prime }+r^{\prime })+1$
and $\mathcal{E}$ is bounded as above. Therefore 
\begin{equation*}
\Delta =(p+r)+(q+s)-(p^{\prime }+r^{\prime })-8
\end{equation*}%
where $p+r=p_{1}(p^{\prime }+r^{\prime })$, $q+s=q_{1}(p^{\prime }+r^{\prime
})$ and $2\leq p_{1}<q_{1}$. We find that the only possibility for $\Delta
\leq 0$ is $p+r=4$, $q+s=6$. But then Lemma 2.2 (c) gives $\sum \eta _{i}>3$%
. Repeating the above procedure, we get $\mathcal{E}\leq 12$ and $2\delta
_{\max }=14$. \medskip

\subsection{Case A2.}

Here we assume that the contribution from finite singular points is small
and the contribution from infinity is maximal possible, so that $ext\nu
_{1}=\dots =ext\nu _{4}=1$ and $\nu _{\infty }$ is maximal (the case with $%
\nu _{0}$ maximal is analogous). By Proposition 2.6 (inequality (2.9)) and
Lemma 2.2 (a) we get 
\begin{equation*}
\nu _{\infty }\leq p+r+q+s-7.
\end{equation*}%
Hence, by formula (2.6) for $ps\not=qr,$ 
\begin{equation*}
\mathcal{E}\leq p^{\prime }(p+r+q+s-7)+8
\end{equation*}%
where $p^{\prime }=\gcd (p,q)$ (see (2.5)). Therefore 
\begin{equation*}
\Delta =2\delta _{\max }-\mathcal{E}\geq (p+r-1-p^{\prime })(q+s-1-p^{\prime
})-\left( p^{\prime }\right) ^{2}+5p^{\prime }-8.
\end{equation*}%
We have $p^{\prime }\geq 2$, since otherwise ($p^{\prime }=1$) the
singularity at infinity is quasi--homogeneous and $\nu _{\infty }=0$ by
definition. Therefore it is enough to prove that $(p+r-1-p^{\prime
})(q+s-1-p^{\prime })\geq \left( p^{\prime }\right) ^{2}-1$. This is
obviously true if $p\geq 2p^{\prime }$, since $q\geq p^{\prime }+p$.
Otherwise, the handsomeness property ensures that $r\geq p^{\prime }$, so $%
p+r-1-p^{\prime }\geq p^{\prime }-1$ and $q+s-1-p^{\prime }\geq 2p^{\prime
}-1$. Hence we ask whether $(p^{\prime }-1)(2p^{\prime }-1)\geq \left(
p^{\prime }\right) ^{2}-1$. But this is always true for $p^{\prime }\geq 2$.
\medskip

\subsection{Case B2.}

Let us denote $p+r=e$, $q+s=f$, $\gcd (e,f)=p^{\prime }+r^{\prime
}=e^{\prime }$. From inequality (2.9) in Proposition 2.6 we get (as in case
A2) 
\begin{equation*}
\nu _{\inf }=\nu _{0}+\nu _{\infty }+\nu _{\tan }\leq e+f-7.
\end{equation*}%
Hence $\mathcal{E}\leq e^{\prime }(e+f-6)+8$ (see (2.6)), whereas $2\delta
_{\max }=(e-1)(f-1)+1-e^{\prime }$ (see Proposition 2.4). Therefore $\Delta
\geq (e-e^{\prime }-1)(f-e^{\prime }-1)-7-{e^{\prime }}^{2}+3e^{\prime }$.
Since $e\geq 2e^{\prime }$ and $f\geq 3e^{\prime }$, we get $\Delta >0$ if $%
(e^{\prime })^{2}\geq 7$. So we assume that $e^{\prime }=2$ and then $\Delta
\geq (e-3)(f-3)-5$. If $e\geq 6$ then $f\geq 8$ and we get $\Delta >0$.
Hence $e=4$ (it must be even). But then Lemma 2.2 (c) implies that $\nu
_{\inf }\leq e+f-8$, so $\Delta \geq (e-e^{\prime }-1)(f-e^{\prime }-1)-7-{%
e^{\prime }}^{2}+4e^{\prime }$. We observe that $\Delta >0$, unless $e=4$
and $f=6$ in which case we obtain $\Delta =0$.

We have to exclude the latter possibility. This can be done by computing the
sum of $\delta $--invariants of singularities (numbers of double points) of
the curve $C$ explicitly. If $e=4$ and $f=6$ then $p=r=2$, $q=s=3$. As $%
t\rightarrow \infty $ (respectively $t\rightarrow 0$) we have $x\sim t^{2}$, 
$y\sim t^{3}$ (respectively $x\sim t^{-2}$, $y\sim t^{-3}$). In the local
coordinates $u=x/y$, $w=1/y$, $s=t^{-1}$ and $s\rightarrow 0$ we have $%
u=s+\dots $, $w=s^{3}+\dots $. Thus the both branches are smooth at
infinity. Then $\nu _{0}=\nu _{\infty }=0$ and $\nu _{\inf }=\nu _{\tan }$.
The requirement $\Delta =0$ implies $\nu _{\tan }=2$. Therefore, if we
consider Puiseux expansions 
\begin{align*}
y& =c_{0}x^{3/2}+c_{1}x+c_{2}x^{1/2}+\dots & & \text{as $t\rightarrow \infty
,$} \\
y& =d_{0}x^{3/2}+d_{1}x+d_{2}x^{1/2}+\dots & & \text{as $t\rightarrow 0$},
\end{align*}%
then we must have $c_{0}=d_{0}$ and $c_{1}=d_{1}$ and the codimension bound
prohibits that $c_{2}=d_{2}$ if earlier terms agree. It follows that the
intersection index of the two branches at infinity is $5$. So the $\delta $%
--invariant of the singularity at infinity is $5$. Adding $4$ from the cusps
at finite distance we obtain $9$. But $C$ is rational of degree $6$, so the
sum of its $\delta $--invariants is $\frac{1}{2}5\cdot 4$. Hence it must
have an additional double point at finite distance. \medskip

\subsection{Case C1.}

Similarly as in case A1 (using the bound (2.10) and Lemma 2.2 (a)) we get 
\begin{equation*}
ext\nu _{1}\leq p+q+r+s-5.
\end{equation*}%
We get then $\mathcal{E}\leq (p+r-2)(q+s-1)+6$ and hence $\Delta \geq
(q+s-1)+\det \,^{\prime }-6,$ where $\det \,^{\prime }$\ is defined in Lemma
2.3.

\medskip

\textbf{Lemma 3.1.} \textit{We have} $\det \,^{\prime }\geq p^{\prime
}r^{\prime }+1\geq 2$.

\medskip

\textit{Proof}. As $p\geq q+p^{\prime }$ and $s\geq r+r^{\prime }$, we infer
that $\det \,^{\prime }\geq p^{\prime }r^{\prime }+(q-1)r^{\prime
}+(r-1)p^{\prime }+1\geq p^{\prime }r^{\prime }+1$. $\qed$ \medskip

By Lemma 3.1 we get $\Delta \geq q+s-5$. So, if $q+s\geq 6,$ then we are
done. Suppose $q+s\leq 5$. Then $p+r\leq 5$ and we apply Lemma 2.2 (c), so
that $ext\nu _{1}\leq p+r+q+s-6$ and $\Delta \geq p+r+q+s+\det \,^{\prime
}-9\geq 1$. \medskip

\subsection{Case C2.}

As in case B2 we assume that $\nu _{\infty }$ is maximal. We get $\nu
_{\infty }\leq p+r+q+s-6$ so $\mathcal{E}\leq p^{\prime }(p+r+q+s-6)+8$ and 
\begin{equation*}
\Delta \geq (p+r-1-p^{\prime })(q+s-1-p^{\prime })-\left( p^{\prime }\right)
^{2}+4p^{\prime }-8+\det \,^{\prime }.
\end{equation*}%
As $p^{\prime }\geq 2$, using Lemma 3.1, it suffices to show that $%
(p+r-p^{\prime }-1)(q+s-p^{\prime }-1)\geq \left( p^{\prime }\right)
^{2}-p^{\prime }$. If $q\geq 2p^{\prime }$ then $p\geq 3p^{\prime }$. Thus $%
p+r-p^{\prime }-1\geq 2p^{\prime }$ and $q+s-p^{\prime }-1\geq p^{\prime }+1$
($s>r>0$) and we are done. So assume that $q=p^{\prime }$. By the
handsomeness $s\geq p^{\prime }$ (see Definition 2.2 and Lemma 2.4). This
implies that $p+r-p^{\prime }-1\geq p^{\prime }$ and $q+s-p^{\prime }-1\geq
p^{\prime }-1$. Hence $\Delta >0.$

\medskip

\subsection{Case D1.}

We use the bound $\sum ext\nu _{i}+\nu _{\infty }\leq p-|r|+q+s+2+\lfloor
\left( \left\vert r\right\vert -1\right) /s\rfloor -\sum \eta _{i}$ from
(2.11). We will treat only the case $p+r\leq q+s$; the computations are
almost identical in the opposite case.

Assume firstly that $m_{1}\geq m_{2}\geq 3$. Then $m_{1}\leq p+r-3$ and $%
p+r\geq 6$. Moreover, $ext\nu _{2}\geq 3$ (because the coefficients before $%
(t-t_{2})$ and $(t-t_{2})^{2}$ in both $\varphi (t)$ and $\psi (t)$ in (2.3)
must vanish). Therefore $ext\nu _{1}\leq p-|r|+q+s-6+\lfloor \left(
\left\vert r\right\vert -1\right) /s\rfloor $. Hence 
\begin{equation*}
\mathcal{E}=(p-|r|-3)\left( q+s-1+\intfrac{|r|-1}{s}\right)
+3\cdot (3-1)+2+2.
\end{equation*}%
Thus 
\begin{equation*}
\Delta =2(q+s)-p^{\prime }-r^{\prime }-11+ps+|r|q-(p-|r|-3)\cdot
\intfrac{|r|-1}{s}.
\end{equation*}%
Now, since $p-\left\vert r\right\vert \leq q+s,$ we get 
\begin{equation}
|r|q-(p-|r|)\cdot \intfrac{|r|-1}{s}\geq 1-|r|+|r|q\left( 1-%
\frac{1}{s}\right) +\frac{q}{s}.  \label{3.1}
\end{equation}%
Substituting this into $\Delta $ we obtain 
\begin{equation*}
\Delta \geq 2(q+s)-p^{\prime }-r^{\prime }-10+ps-|r|+|r|q\left( 1-\frac{1}{s}%
\right) +\frac{q}{s}.
\end{equation*}%
But $|r|q(1-\frac{1}{s})+\frac{q}{s}\geq 0$. It follows that 
\begin{equation*}
\Delta \geq (q+s)+(p-|r|)-10+p(s-1)+(q+s-p^{\prime }-r^{\prime })\geq 2,
\end{equation*}%
since $q+s\geq 6$, $p-|r|\geq 6$ and the last two terms in the above formula
are non-negative.

We are left with the case $m_{2}=2$. Then $m_{3}=m_{4}=2$, so $\sum \eta
_{i}>3$ by Lemma 2.2 (b). We obtain $ext\nu _{1}\leq p-|r|+q+s-5+\lfloor
\left( \left\vert r\right\vert -1\right) /s\rfloor $. Hence 
\begin{equation*}
\mathcal{E}=(p-|r|-2)\left( q+s-1+\intfrac{|r|-1}{s}\right)
+6.
\end{equation*}%
Therefore 
\begin{equation*}
\Delta =q+s-p^{\prime }-r^{\prime }-6+ps+|r|q-(p-|r|-2)\cdot \intfrac{|r|-1}{s}.
\end{equation*}%
If $s=1$ then $\Delta =q-p^{\prime }-8+|r|(q-p+|r|+3)+2p\geq q+2p-p^{\prime
}+2\left\vert r\right\vert -8$, because $p-|r|\leq q+s=q+1$. But $p^{\prime
}\leq p/2$ and $p\geq 5$, as $p-|r|\geq 4$. So $\Delta >0$.

Finally, let us assume that $s\geq 2$. By (3.1) we have 
\begin{equation*}
\Delta \geq q+s-p^{\prime }-r^{\prime }-5+ps-|r|.
\end{equation*}%
Then $q+s-p^{\prime }-r^{\prime }\geq 0$ and $ps-|r|\geq p+p-|r|\geq 9$, so $%
\Delta \geq 4$. \medskip

\subsection{Case D2.}

Here $\nu _{\infty }$ (or $\nu _{0}$) is bounded from above by $%
p-|r|+q+s+\lfloor \left( \left\vert r\right\vert -1\right) /s\rfloor +2-7$.
Assume that $p^{\prime }\geq r^{\prime }$. It follows that $\mathcal{E}\leq
8+p^{\prime }(p-|r|+q+s+\lfloor \left( \left\vert r\right\vert -1\right)
/s\rfloor -5)$. So

\begin{equation*}
\begin{split}
\Delta =& (p-|r|-1)(q+s-1)+ps+|r|q-p^{\prime }-r^{\prime }-7 \\
& -p^{\prime }\left( p-|r|+q+s+\intfrac{|r|-1}{s}-5\right)
\end{split}%
\end{equation*}%
This can be transformed into

\begin{equation}
\begin{split}
\Delta =& (p-p^{\prime }-1)(q+s-p^{\prime }-1)+(p-|r|)s \\
& +\left( |r|-p^{\prime }-\intfrac{|r|-1}{s}+2\right)
p^{\prime }+(|r|-r^{\prime })-7.
\end{split}
\label{3.2}
\end{equation}

Assume that $q\geq 2p^{\prime }$. Obviously, also $p\geq 2p^{\prime }$ and $%
p\neq q$; thus either $p$ or $q$ is at least $3p^{\prime }$. Then $%
(p-p^{\prime }-1)(q+s-p^{\prime }-1)\geq 2p^{\prime }(p^{\prime }-1)$.
Moreover, $(p-|r|)s\geq 4$ and $|r|-p^{\prime }-\lfloor \left( \left\vert
r\right\vert -1\right) /s\rfloor +2\geq 3-p^{\prime }$. Hence $\Delta \geq
\left( p^{\prime }\right) ^{2}+p^{\prime }-3>0$, as $p^{\prime }\geq 2$.

Therefore $q=p^{\prime }$. By the handsomeness $s\geq p^{\prime }$ (see
Definition 2.2 and Lemma 2.4). It follows that $(p-p^{\prime
}-1)(q+s-p^{\prime }-1)\geq (p^{\prime }-1)^{2}$. Thus 
\begin{equation*}
\Delta \geq (p^{\prime }-1)^{2}+(p-|r|)s+(3-p^{\prime })p^{\prime }-7.
\end{equation*}%
But $(p-|r|)s\geq 8$, so $\Delta >0$.

Now let us turn to the case $r^{\prime }>p^{\prime }$. Equation (3.2) then
becomes

\begin{equation*}
\begin{split}
\Delta \geq & (p-r^{\prime }-1)(q+s-r^{\prime }-1)+(p-|r|)s \\
& +\left( |r|-r^{\prime }-\intfrac{|r|-1}{s}+2\right)
r^{\prime }+|r|-p^{\prime }-7.
\end{split}%
\end{equation*}%
Here $s\geq r^{\prime }\geq 2$. If $|r|\geq 2r^{\prime }$ then $%
|r|-r^{\prime }-\lfloor \left( \left\vert r\right\vert -1\right) /s\rfloor
>0 $. Since $p-r^{\prime }\geq p-|r|$ we have $(p-r^{\prime
}-1)(q+s-r^{\prime }-1)\geq 0$. We infer that $\Delta \geq
(p-|r|)s+2r^{\prime }+|r|-p^{\prime }-7>1$ for $|r|>r^{\prime }>p^{\prime }$%
. So let $|r|=r^{\prime }$. Then $|r|\leq s$, so $\lfloor \left( \left\vert
r\right\vert -1\right) /s\rfloor =0 $. Hence $\Delta \geq (p-|r|)s-7>0$.

\medskip

\subsection{Case E1.}

We have here 
\begin{equation*}
\sum ext\nu _{i}+\nu _{\infty }\leq p-|r|+s-|q|+\lfloor \left( \left\vert
q\right\vert -1\right) /p\rfloor +\lfloor \left( \left\vert s\right\vert
-1\right) /r\rfloor +3-\sum \eta _{i},
\end{equation*}%
i.e. the bound (2.12) holds. It is easy to observe that at most one of the $%
\lfloor \left( \left\vert q\right\vert -1\right) /p\rfloor $ and $\lfloor
\left( \left\vert s\right\vert -1\right) /r\rfloor $ can be non--zero.
Following [BZ2] we introduce the quantities $K=p-|r|$ and $L=s-|q|$ with $%
K\leq L$.

\medskip

\emph{Subcase (i): }$\lfloor \left( \left\vert q\right\vert -1\right)
/p\rfloor $\emph{$>0$.} Putting $ext\nu _{i}=1$ for $i\geq 2,$ we get $%
ext\nu _{1}<K+L+\lfloor \left( \left\vert q\right\vert -1\right) /p\rfloor
+3-3-2.$ Since all terms in this inequality are integers, we have $ext\nu
_{1}\leq K+L+\lfloor \left( \left\vert q\right\vert -1\right) /p\rfloor -3$
and $m_{1}\leq K-2$. So 
\begin{equation*}
\mathcal{E}\leq (K-2)\left( L+\intfrac{|q|-1}{p}+1\right) +6.
\end{equation*}%
Since $2\delta _{\max }=(K-1)(L-1)+|r|L+|q|K+KL-p^{\prime }-r^{\prime }+1$
we get 
\begin{equation}
\Delta \geq KL+K(|q|-2)+L(|r|+1)-(K-2)\cdot \intfrac{|q|-1}{p}-p^{\prime }-r^{\prime }-2.  \label{3.3}
\end{equation}%
Now $p>K-2$. Therefore $(K-2)\cdot \lfloor \left( \left\vert q\right\vert
-1\right) /p\rfloor \leq |q|-1$. Thus 
\begin{equation*}
\Delta \geq KL-3+(K-1)(|q|-2)+L(|r|+1)-p^{\prime }-r^{\prime }.
\end{equation*}%
The above inequality can be rewritten as 
\begin{equation*}
\Delta \geq KL-4+(K-2)(|q|-2)+(L-1)(|r|+1)+(|q|-p^{\prime })+(|r|-r^{\prime
}).
\end{equation*}%
We have $K,L\geq 4$ and $|q|\geq 3$. Therefore 
\begin{equation*}
\Delta \geq KL-4-(K-2)+3\cdot 2>0.
\end{equation*}%
\medskip

\emph{Subcase (ii): }$\lfloor \left( \left\vert q\right\vert -1\right)
/p\rfloor $\emph{$=0$.} Then we get an equation similar to (3.3) 
\begin{equation*}
\Delta \geq KL+K(|q|-2)+L(|r|+1)-(K-2)\cdot \intfrac{|r|-1}{s}-p^{\prime }-r^{\prime }-2.
\end{equation*}%
We have $K-2\leq L-2<s$ and hence $(K-2)\cdot \lfloor \left( \left\vert
r\right\vert -1\right) /s\rfloor \leq |r|-1$. Using this we transform the
above inequality into 
\begin{equation*}
\Delta \geq KL-1+(K-1)(|q|-2)+(L-2)(|r|+1)+(|q|-p^{\prime })+(|r|-r^{\prime
}).
\end{equation*}%
As $|q|\geq 1$ we get $\Delta >0$. \medskip

\subsection{Case E2.}

Assume that $p^{\prime }\geq r^{\prime }$. We will not impose, however, the
inequality $K\leq L$. Then $\mathcal{E}\leq p^{\prime }(K+L-4+\lfloor \left(
\left\vert q\right\vert -1\right) /p\rfloor +\lfloor \left( \left\vert
r\right\vert -1\right) /s\rfloor )+8$. On the other hand, $2\delta _{\max
}=KL-K-L+2+K|q|+pL-p^{\prime }-r^{\prime }$. Henceforth 
\begin{equation*}
\Delta \geq KL-K-L-6+(p-p^{\prime })L+(|q|-p^{\prime })K-p^{\prime }\left(
\intfrac{|q|-1}{p}+\intfrac{|r|-1}{s}-2\right) .
\end{equation*}%
If $|r|-1<s$ then $p^{\prime }\cdot \lfloor \left( \left\vert q\right\vert
-1\right) /p\rfloor \leq |q|-1$. Hence we are left with $\Delta \geq
KL-K-L-5+p^{\prime }+(|q|-p^{\prime })(K-1),$ where the latter expression is
positive.

Therefore $|r|-1\geq s=|q|+L$. Since $p^{\prime }\leq |q|<s$ we infer that $%
p^{\prime }\cdot \lfloor \left( \left\vert r\right\vert -1\right) /s\rfloor
\leq |r|-1$. Reminding that $(p-p^{\prime })L>(p-p^{\prime })$ we obtain 
\begin{equation*}
\Delta \geq KL-L-5+(|q|-p^{\prime })K+p^{\prime }+(p-K-|r|).
\end{equation*}%
As $K+|r|=p$ we get $\Delta \geq (K-1)L-5+(|q|-p^{\prime })K+p^{\prime }$.
Since $(K-1)L\geq 12$ we get $\Delta >0$.\medskip

Now the proof of Main Theorem is complete. $\qed$

\end{document}